\documentclass[12pt]{article} 

\usepackage{latexsym, 
amscd, 
graphicx, amssymb}

\newtheorem{thm}{Theorem}[section]

\newtheorem{prop}[thm]{Proposition}

\newtheorem{rema}[thm]{Remark}

\newcommand{\halmos}{\rule{1ex}{1.4ex}}

\newcommand{\text}[1]{\mbox{\rm #1}}

\newcommand{\nn}{\nonumber \\}

 \newcommand{\epfv}{\hspace*{\fill}\mbox{$\halmos$}\vspace{1em}}

\newcommand{\tr}{\mbox{\rm Tr}}
\newcommand{\A}{\mathcal{A}}
\newcommand{\Y}{\mathcal{Y}}
\newcommand{\C}{\mathbb{C}}
\newcommand{\Z}{\mathbb{Z}}

\newcommand{\Q}{\mathbb{Q}}
\newcommand{\N}{\mathbb{N}}

\title{ {\bf Vertex operator algebras, the Verlinde conjecture and 
modular tensor categories} }
\date{}
\author{Yi-Zhi Huang}

\begin{document}

\bibliographystyle{alpha}
\maketitle

\begin{abstract} 
Let $V$ be a simple vertex operator algebra satisfying the 
following conditions: (i) $V_{(n)}=0$ for $n<0$, 
$V_{(0)}=\mathbb{C}\mathbf{1}$ and the contragredient module
$V'$ is isomorphic to $V$ as a $V$-module. (ii)
Every $\mathbb{N}$-gradable weak $V$-module is completely 
reducible. (iii) $V$ is $C_{2}$-cofinite. 
We announce a proof of the Verlinde 
conjecture for $V$, that is, of the statement that 
the matrices formed by the fusion rules 
among irreducible $V$-modules are diagonalized by 
the matrix given by the action of the modular transformation
$\tau\mapsto -1/\tau$ on the space of  characters
of irreducible $V$-modules. We discuss some 
consequences of the Verlinde conjecture, including
the Verlinde formula for the fusion rules, a formula 
for the matrix given by the action of $\tau\mapsto -1/\tau$
and the symmetry of this
matrix. We also announce 
a proof of the rigidity and nondegeneracy property of the 
braided tensor category structure on the category
of $V$-modules when $V$ satisfies in addition the condition that 
irreducible $V$-modules not equivalent to $V$ has no nonzero 
elements of weight $0$.
In particular, the category of $V$-modules 
has a natural structure of modular tensor category. 
\end{abstract}

\renewcommand{\theequation}{\thesection.\arabic{equation}}
\renewcommand{\thethm}{\thesection.\arabic{thm}}
\setcounter{equation}{0}
\setcounter{thm}{0}
\setcounter{section}{-1}

\section{Introduction}

In 1987, by comparing fusion algebras with certain algebras obtained in
the study of conformal field theories on genus-one Riemann surfaces,
Verlinde \cite{V} conjectured that the matrices formed by the fusion
rules are diagonalized by the matrix given
by the action of the modular transformation $\tau\mapsto -1/\tau$ on the
space of characters of a rational conformal field theory.
From this conjecture,
Verlinde obtained the famous Verlinde formulas for the fusion rules and,
more generally, for the dimensions of conformal blocks on Riemann
surfaces of arbitrary genera.  In the particular case of the conformal
field theories associated to affine Lie algebras (the
Wess-Zumino-Novikov-Witten models), the Verlinde formulas give a
surprising formula for the dimensions of the spaces of sections of the
``generalized theta divisors"; this has given rise to a great deal of
excitement and new mathematics. See the works \cite{TUY} by
Tsuchiya-Ueno-Yamada, \cite{BL} by Beauville-Laszlo, \cite{F} by
Faltings and \cite{KNR} by Kumar-Narasimhan-Ramanathan for details and
proofs of this particular case of the Verlinde formulas. 

In 1988, Moore and Seiberg \cite{MS1} showed on a physical level of
rigor that the Verlinde conjecture is a consequence of the axioms for
rational conformal field theories. This result of Moore and Seiberg is
based on certain polynomial equations which they derived from the axioms for
rational conformal field theories \cite{MS1} \cite{MS2}. Moore and
Seiberg further demonstrated that these polynomial equations are
actually conformal-field-theoretic analogues of the tensor category
theory for group representations.  This work of Moore and Seiberg
greatly advanced our understanding of the structure of conformal field
theories. In particular, the notion of modular tensor category was later
introduced to summarize the properties of the Moore-Seiberg polynomial
equations and has played a central role in the developments of
conformal field theories and three-dimensional topological field
theories.  See for example \cite{T} and \cite{BK} for the theory of
modular tensor categories, their applications and references to many
important works done by mathematicians and physicists. 

The work of Moore and Seiberg gave a conceptual 
understanding of the Verlinde conjecture and the modular tensor 
categories arising in conformal field theories. However,
it is a very hard problem to mathematically
construct theories satisfying the axioms for rational 
conformal field theories. In fact, these axioms for rational
conformal field theories are much stronger than the Verlinde 
conjecture and the modular tensor category structures.
In the general theory of vertex operator algebras, introduced and 
studied first by Borcherds \cite{B} and Frenkel-Lepowsky-Meurman
\cite{FLM}, a mathematical version of the notion of fusion 
rule was introduced and studied by Frenkel, Lepowsky and the author 
in \cite{FHL} using intertwining operators,
and the modular transformations were given by Zhu's modular 
invariance theorem \cite{Z}. Using these notions and some natural 
conditions, including in particular 
Zhu's $C_{2}$-cofiniteness condition, one can formulate
a general version of the Verlinde conjecture
in the framework of the theory of vertex operator algebras.
Further results on intertwining operators and modular invariance
were obtained in \cite{HL1}--\cite{HL4} by Huang-Lepowsky, in
\cite{H1}, \cite{H2} and \cite{H4} by the author, in \cite{DLM} by
Dong-Li-Mason and in \cite{M} by Miyamoto. But these results were still
not enough for the proof of this general version of the Verlinde
conjecture.  The main obstructions were the duality and modular
invariance properties for genus-zero and genus-one {\it multi-point}
correlation functions constructed from intertwining operators for a
vertex operator algebra satisfying the conditions above. These
properties have recently been proved in \cite{H5} and \cite{H6}.

In this paper, we announce a proof of the general version of 
the Verlinde conjecture above. 
Our theorem assumes
only that the vertex operator algebra that we consider satisfies certain
natural grading, finiteness and reductivity properties (see Section 2).
We also discuss some  consequences of our theorem, including
the Verlinde formula for the fusion rules, a formula 
for the matrix given by the action of $\tau\mapsto -1/\tau$
and the symmetry of this
matrix. For the details, see \cite{H7}.
We also announce a proof of the rigidity and nondegenracy condition
of the braided tensor category structure on the category of 
modules for such a vertex operator algebra constructed by Lepowsky 
and the author \cite{HL1}--\cite{HL4} \cite{H1} \cite{H5}, when 
$V$ satisfies in addition the condition that 
irreducible modules not equivalent to the algebra (as a module)
has no nonzero 
elements of weight $0$. 
In particular, the category of modules for such a vertex 
operator algebra has a natural structure of modular tensor category.

This paper is organized as follows: In Section 1, we give the
definitions of fusion rule, of the fusing and of the braiding
isomorphisms in terms of matrix elements, and 
of the corresponding action of the modular transformation.  
These are the basic
ingredients needed in the formulations of the main results given in
Sections 2 and 3 and they are in fact based on
substantial mathematical results in \cite{H1}, 
\cite{H2}, \cite{H4}, and in \cite{Z}
and \cite{DLM}, respectively. Our main theorems on the Verlinde
conjecture, on the Verlinde formula for the fusion rules, on the formula for
the matrix given by the action of $\tau\mapsto -1/\tau$, and on the symmetry
of this matrix, are stated in Section 2. A very brief sketch of the 
proof of the Verlinde conjecture is given in this section. 
In Section 3, our main theorem on the modular tensor category structure
is stated and a sketch of the proof is given.

\paragraph{Acknowledgment} I am grateful to Jim Lepowsky and Robert Wilson 
for comments. The author is partially supported 
by NSF grant DMS-0401302.

\renewcommand{\theequation}{\thesection.\arabic{equation}}
\renewcommand{\thethm}{\thesection.\arabic{thm}}
\setcounter{equation}{0}
\setcounter{thm}{0}

\section{Fusion rules, fusing and braiding isomorphisms
and modular transformations}

We assume that the reader is familiar with the basic definitions and 
results in the theory of vertex operator algebras as introduced and 
presented in \cite{B} and \cite{FLM}. We shall use the notations, 
terminology and formulations in 
\cite{FLM}, \cite{FHL} and \cite{LL}.

Let $V$ be a simple vertex operator algebra, $V'$ the contragredient 
module of $V$, and $C_{2}(V)$ the subspace
of $V$ spanned by $u_{-2}v$ for $u, v\in V$. In the present paper, we
shall always assume that $V$ satisfies the following conditions:

\begin{enumerate}

\item $V_{(n)}=0$ for $n<0$, $V_{(0)}=\mathbb{C}\mathbf{1}$ and 
$V'$ is isomorphic to $V$ as a $V$-module.

\item Every $\mathbb{N}$-gradable weak $V$-module is completely 
reducible.

\item $V$ is $C_{2}$-cofinite, that is, $\dim V/C_{2}(V)<\infty$. 

\end{enumerate}

We recall that an $\mathbb{N}$-gradable 
weak $V$-module 
is a 
vector space that admits  an $\mathbb{N}$-grading
$W=\coprod_{n\in \mathbb{N}}W_{[n]}$, equipped with 
a vertex operator map
\begin{eqnarray*}
Y: V\otimes W&\to &W[[z, z^{-1}]] \\
u\otimes w&\mapsto &Y(u, z)w=\sum_{n\in \mathbb{Z}}u_{n}z^{-n-1}
\end{eqnarray*}
satisfying all axioms for $V$-modules except that the condition
$L(0)w=nw$ for $w\in 
W_{(n)}$ is replaced by $u_{k}w\in W_{[m-k-1+n]}$ for 
$u\in V_{(m)}$ and $w\in W_{[n]}$. Condition 2 is equivalent to
the statement that 
every finitely-generated $\mathbb{N}$-gradable weak $V$-module
is a $V$-module and every $V$-module is completely reducible.

From \cite{DLM}, we know that there are only finitely many
inequivalent irreducible $V$-modules.  Let $\mathcal{A}$ be the set of
equivalence classes of irreducible $V$-modules. We denote the
equivalence class containing $V$ by $e$. For each $a\in \mathcal{A}$,
we choose a representative $W^{a}$ of $a$. Note that the contragredient
module of an irreducible module is also irreducible (see \cite{FHL}).
So we have a map
\begin{eqnarray*}
^{\prime}:
\mathcal{A}&\to& \mathcal{A}\\
a&\mapsto& a'.
\end{eqnarray*}
From \cite{AM} and \cite{DLM}, we know that irreducible $V$-modules 
are in fact graded by rational numbers. 
Thus for $a\in \A$, 
there exist $h_{a}\in \Q$ such that 
$W^{a}=\coprod_{n\in h_{a}+\N}W^{a}_{(n)}$.

Let $\mathcal{V}_{a_{1}a_{2}}^{a_{3}}$ for $a_{1}, a_{2}, a_{3}
\in \mathcal{A}$
be the space of intertwining operators of 
type ${W^{a_{3}}\choose W^{a_{1}}W^{a_{2}}}$ and 
$N_{a_{1}a_{2}}^{a_{3}}$ for $a_{1}, a_{2}, a_{3}
\in \mathcal{A}$ the fusion rule, that is, the dimension of 
the space of intertwining operators of 
type ${W^{a_{3}}\choose W^{a_{1}}W^{a_{2}}}$. For any 
$\Y\in \mathcal{V}_{a_{1}a_{2}}^{a_{3}}$, we know from \cite{FHL}
that 
for $w_{a_{1}}\in W^{a_{1}}$ and $w_{a_{2}}\in W^{a_{2}}$
\begin{equation}\label{y-delta}
\Y(w_{a_{1}}, x)w_{a_{2}}\in x^{\Delta(\Y)}W^{a_{3}}[[x, x^{-1}]],
\end{equation}
where 
$$\Delta(\mathcal{Y})=h_{a_{3}}-h_{a_{1}}-h_{a_{2}}.$$

From \cite{GN},
\cite{L},
\cite{AN},
\cite{H6}, we also know that the fusion rules 
$N_{a_{1}a_{2}}^{a_{3}}$ for 
$a_{1}, a_{2}, a_{3}\in \mathcal{A}$ are all finite. 
For $a\in \A$, let $\mathcal{N}(a)$ be the matrix 
whose entries are $N_{aa_{1}}^{a_{2}}$ for $a_{1}, a_{2}\in \A$,
that is, 
$$\mathcal{N}(a)=(N_{aa_{1}}^{a_{2}}).$$

We also need matrix elements of fusing and braiding isomorphisms. 
In the proof of the Verlinde conjecture, we need 
to use several bases of one space of intertwining 
operators. We shall  use $p=1, 2, 3, 4, 5, 6, \dots$
to label different bases. 
For $p=1, 2, 3, 4, 5, 6, \dots$ and $a_{1}, a_{2}, 
a_{3}\in \mathcal{A}$, let 
$\{\Y_{a_{1}a_{2}; i}^{a_{3}; (p)}\;|\; i=1, \dots, 
N_{a_{1}a_{2}}^{a_{3}}\}$, be a
basis of $\mathcal{V}_{a_{1}a_{2}}^{a_{3}}$. 
For $a_{1}, \dots, a_{6}\in \mathcal{A}$, 
$w_{a_{1}}\in W^{a_{1}}$, $w_{a_{2}}\in W^{a_{2}}$,
$w_{a_{3}}\in W^{a_{3}}$, and $w'_{a_{4}}\in (W^{a_{4}})'$,
using the differential equations satisfied by the series
$$\langle w'_{a_{4}}, \Y_{a_{1}a_{5}; i}^{a_{4}; (1)}(w_{a_{1}}, x_{1})
\Y_{a_{2}a_{3}; j}^{a_{5}; (2)}(w_{a_{2}}, x_{2})w_{a_{3}}\rangle
|_{x_{1}^{n}=e^{n\log z_{1}}, \; x_{2}^{n}=e^{n\log z_{2}}, \;
n\in \Q}$$
and 
$$\langle w'_{a_{4}}, 
\Y_{a_{6}a_{3}; k}^{a_{4}; (3)}
(\Y_{a_{1}a_{2}; l}^{a_{6}; (4)}(w_{a_{1}}, x_{0})
w_{a_{2}}, x_{2})w_{a_{3}}\rangle
|_{x_{0}^{n}=e^{n\log (z_{1}-z_{2})}, \; x_{2}^{n}=e^{n\log z_{2}}, \;
n\in \Q},$$
it was proved in \cite{H5} that these series are convergent 
in the regions $|z_{1}|>|z_{2}|>0$ and $|z_{2}|>|z_{1}-z_{2}|>0$,
respectively. Note that for any $a_{1}$, $a_{2}$, $a_{3}$, $a_{4}$, $a_{5}$,
$a_{6}\in \A$, 
$\{\Y_{a_{1}a_{5}; i}^{a_{4}; (1)}
\otimes \Y_{a_{2}a_{3}; j}^{a_{5}; (2)}\;|\;i=1, \dots, 
N_{a_{1}a_{5}}^{a_{4}}, j=1, \dots, N_{a_{2}a_{3}}^{a_{5}}\}$
and 
$\{\Y_{a_{6}a_{3}; l}^{a_{4}; (3)}
\otimes \Y_{a_{1}a_{2}; k}^{a_{6}; (4)}\;|\; l=1, \dots, 
N_{a_{6}a_{3}}^{a_{4}}, k=1, \dots, N_{a_{1}a_{2}}^{a_{6}}\}$
are bases of $\mathcal{V}_{a_{1}a_{5}}^{a_{4}}\otimes 
\mathcal{V}_{a_{2}a_{3}}^{a_{5}}$ and 
$\mathcal{V}_{a_{6}a_{3}}^{a_{4}}
\otimes \mathcal{V}_{a_{1}a_{2}}^{a_{6}}$, respectively. 
The associativity of intertwining operators proved and 
studied in \cite{H1}, \cite{H4} and
\cite{H5} says that 
there exist 
$$F(\Y_{a_{1}a_{5}; i}^{a_{4}; (1)}\otimes \Y_{a_{2}a_{3}; j}^{a_{5}; (2)}; 
\Y_{a_{6}a_{3}; l}^{a_{4}; (3)}
\otimes \Y_{a_{1}a_{2}; k}^{a_{6}; (4)}) \in \C$$
for $a_{1}, \dots, a_{6}\in \mathcal{A}$, $i=1, \dots, 
N_{a_{1}a_{5}}^{a_{4}}$, $j=1, \dots, 
N_{a_{2}a_{3}}^{a_{5}}$, $k=1, \dots, 
N_{a_{6}a_{3}}^{a_{4}}$, $l=1, \dots, 
N_{a_{1}a_{2}}^{a_{6}}$
such that 
\begin{eqnarray}\label{assoc}
\lefteqn{\langle w'_{a_{4}}, \Y_{a_{1}a_{5}; i}^{a_{4}; (1)}(w_{a_{1}}, x_{1})
\Y_{a_{2}a_{3}; j}^{a_{5}; (2)}(w_{a_{2}}, z_{2})w_{a_{3}}\rangle
|_{x_{1}^{n}=e^{n\log z_{1}}, \; x_{2}^{n}=e^{n\log z_{2}}, \;
n\in \Q}}\nn
&&=\sum_{a_{6}\in \A}
\sum_{k=1}^{N_{a_{6}a_{3}}^{a_{4}}}\sum_{l=1}^{N_{a_{1}a_{2}}^{a_{6}}}
F(\Y_{a_{1}a_{5}; i}^{a_{4}; (1)}\otimes \Y_{a_{2}a_{3}; j}^{a_{5}; (2)}; 
\Y_{a_{6}a_{3}; l}^{a_{4}; (3)}
\otimes \Y_{a_{1}a_{2}; k}^{a_{6}; (4)})\cdot\nn
&&\quad \cdot 
\langle w'_{a_{4}}, 
\Y_{a_{6}a_{3}; k}^{a_{4}; (3)}(\Y_{a_{1}a_{2}; l}^{a_{6}; (4)}(w_{a_{1}}, z_{1}-z_{2})
w_{a_{2}}, z_{2})w_{a_{3}}\rangle
|_{x_{0}^{n}=e^{n\log (z_{1}-z_{2})}, \; x_{2}^{n}=e^{n\log z_{2}}, \;
n\in \Q}\nn
&&
\end{eqnarray}
when $|z_{1}|>|z_{2}|>|z_{1}-z_{2}|>0$, 
for $a_{1}, \dots, a_{5}\in \A$, $w_{a_{1}}\in W^{a_{1}}$, 
$w_{a_{2}}\in W^{a_{2}}$, $w_{a_{3}}\in W^{a_{3}}$,  
$w'_{a_{4}}\in (W^{a_{4}})'$, $i=1, \dots, 
N_{a_{1}a_{5}}^{a_{4}}$ and $j=1, \dots, 
N_{a_{2}a_{3}}^{a_{5}}$. The numbers 
$$F(\Y_{a_{1}a_{5}; i}^{a_{4}; (1)}\otimes \Y_{a_{2}a_{3}; j}^{a_{5}; (2)}; 
\Y_{a_{6}a_{3}; k}^{a_{4}; (3)}
\otimes \Y_{a_{1}a_{2}; l}^{a_{6}; (4)})$$
together give a matrix which represents a linear 
isomorphism
$$\coprod_{a_{1}, a_{2}, a_{3}, a_{4}, a_{5}\in \A}
\mathcal{V}_{a_{1}a_{5}}^{a_{4}}\otimes 
\mathcal{V}_{a_{2}a_{3}}^{a_{5}}\to 
\coprod_{a_{1}, a_{2}, a_{3}, a_{4}, a_{6}\in \A}
\mathcal{V}_{a_{6}a_{3}}^{a_{4}}
\otimes \mathcal{V}_{a_{1}a_{2}}^{a_{6}},$$
called the {\it fusing isomorphism}, 
such that these numbers 
are the matrix elements.

By the commutativity of intertwining operators proved 
and studied in \cite{H2},
\cite{H4}
and \cite{H5}, for any fixed $r\in \Z$, 
there exist 
$$B^{(r)}(\Y_{a_{1}a_{5}; i}^{a_{4}; (1)}
\otimes \Y_{a_{2}a_{3}; j}^{a_{5}; (2)}; 
\Y_{a_{2}a_{6}; l}^{a_{4}; (3)}
\otimes \Y_{a_{1}a_{3}; k}^{a_{6}; (4)}) \in \C$$
for  $a_{1}, \dots, a_{6}\in \mathcal{A}$, $i=1, \dots, 
N_{a_{1}a_{5}}^{a_{4}}$, $j=1, \dots, 
N_{a_{2}a_{3}}^{a_{5}}$, $k=1, \dots, 
N_{a_{2}a_{6}}^{a_{4}}$, $l=1, \dots, 
N_{a_{1}a_{3}}^{a_{6}}$, 
such that 
the analytic extension of the single-valued analytic function
$$\langle w'_{a_{4}}, \Y_{a_{1}a_{5}; i}^{a_{4}; (1)}(w_{a_{1}}, x_{1})
\Y_{a_{2}a_{3}; j}^{a_{5}; (2)}(w_{a_{2}}, x_{2})w_{a_{3}}\rangle
|_{x_{1}^{n}=e^{n\log z_{1}}, \; x_{2}^{n}=e^{n\log z_{2}}, \;
n\in \Q}$$
on the region $|z_{1}|>|z_{2}|>0$, $0\le \arg z_{1}, \arg z_{2}<2\pi$
along the path 
$$t \mapsto \left(\frac{3}{2}
-\frac{e^{(2r+1)\pi i t}}{2}, \frac{3}{2}
+\frac{e^{(2r+1)\pi i t}}{2}\right)$$ 
to the region $|z_{2}|>|z_{1}|>0$, $0\le \arg z_{1}, \arg z_{2}<2\pi$
is 
\begin{eqnarray*}
\lefteqn{\sum_{a_{6}\in \A}
\sum_{k=1}^{N_{a_{2}a_{6}}^{a_{4}}}\sum_{l=1}^{N_{a_{1}a_{3}}^{a_{6}}}
B^{(r)}(\Y_{a_{1}a_{5}; i}^{a_{4}; (1)}\otimes \Y_{a_{2}a_{3}; j}^{a_{5}; (2)};
\Y_{a_{2}a_{6}; k}^{a_{4}; (3)}\otimes 
\Y_{a_{1}a_{3}; l}^{a_{6}; (4)})\cdot}\nn
&&\quad\quad\quad\quad\cdot 
\langle w'_{a_{4}}, 
\Y_{a_{2}a_{6}; k}^{a_{4}; (3)}(
w_{a_{2}}, z_{1})\Y_{a_{1}a_{3}; l}^{a_{6}; (4)}(w_{a_{1}}, z_{2})
w_{a_{3}}\rangle|_{x_{1}^{n}=e^{n\log z_{1}}, \; x_{2}^{n}=e^{n\log z_{2}}, \;
n\in \Q}.
\end{eqnarray*}
The numbers 
$$B^{(r)}(\Y_{a_{1}a_{5}; i}^{a_{4}; (1)}\otimes \Y_{a_{2}a_{3}; j}^{a_{5}; (2)};
\Y_{a_{2}a_{6}; k}^{a_{4}; (3)}\otimes \Y_{a_{1}a_{3}; l}^{a_{6}; (4)})$$
together give a linear isomorphism
$$\coprod_{a_{1}, a_{2}, a_{3}, a_{4}, a_{5}\in \A}
\mathcal{V}_{a_{1}a_{5}}^{a_{4}}\otimes 
\mathcal{V}_{a_{2}a_{3}}^{a_{5}}\to 
\coprod_{a_{1}, a_{2}, a_{3}, a_{4}, a_{6}\in \A}
\mathcal{V}_{a_{2}a_{6}}^{a_{4}}
\otimes \mathcal{V}_{a_{1}a_{3}}^{a_{6}},$$
called the {\it braiding isomorphism}, 
such that these numbers 
are the matrix elements.

We need an action of $S_{3}$ on the space 
$$\mathcal{V}=\coprod_{a_{1}, a_{2}, a_{3}\in 
\mathcal{A}}\mathcal{V}_{a_{1}a_{2}}^{a_{3}}.$$
For $r\in \Z$, $a_{1}, a_{2}, a_{3}\in \mathcal{A}$, consider the
isomorphisms
$\Omega_{r}: \mathcal{V}_{a_{1}a_{2}}^{a_{3}} \to 
\mathcal{V}_{a_{2}a_{1}}^{a_{3}}$ and 
$A_{r}: \mathcal{V}_{a_{1}a_{2}}^{a_{3}} \to 
\mathcal{V}_{a_{1}a'_{3}}^{a'_{2}}$ given in (7.1) and (7.13) 
in \cite{HL2}.
For $a_{1}, a_{2}, a_{3}\in \A$,
$\mathcal{Y}\in \mathcal{V}_{a_{1}a_{2}}^{a_{3}}$, 
we define 
\begin{eqnarray*}
\sigma_{12}(\mathcal{Y})&=&e^{\pi i \Delta(\mathcal{Y})}
\Omega_{-1}(\mathcal{Y})\\
&=&e^{-\pi i \Delta(\mathcal{Y})}\Omega_{0}(\mathcal{Y}),\\
\sigma_{23}(\mathcal{Y})&=&e^{\pi i h_{a_{1}}}
A_{-1}(\mathcal{Y})\\
&=&e^{-\pi i h_{a_{1}}}A_{0}(\mathcal{Y}).
\end{eqnarray*}
We have the following:

\begin{prop}
The actions $\sigma_{12}$ and 
$\sigma_{23}$ of $(12)$ and 
$(23)$ on $\mathcal{V}$
generate a left action of $S_{3}$ on $\mathcal{V}$.
\end{prop} 

We now  choose
a basis $\mathcal{Y}_{a_{1}a_{2}; i}^{a_{3}}$, $i=1, \dots, 
N_{a_{1}a_{2}}^{a_{3}}$, 
of $\mathcal{V}_{a_{1}a_{2}}^{a_{3}}$ for each triple 
$a_{1}, a_{2}, a_{3}\in \mathcal{A}$. 
For $a\in \A$, we choose $\Y_{ea; 1}^{a}$ to be the vertex operator 
$Y_{W^{a}}$ defining the module structure on $W^{a}$ and we choose 
$\Y_{ae; 1}^{a}$ to be the intertwining operator defined using 
the action of $\sigma_{12}$,
\begin{eqnarray*}
\Y_{ae; 1}^{a}(w_{a}, x)u&=&\sigma_{12}(\Y_{ea; 1}^{a})(w_{a}, x)u\nn
&=&e^{xL(-1)}\Y_{ea; 1}^{a}(u, -x)w_{a}\nn
&=&e^{xL(-1)}Y_{W^{a}}(u, -x)w_{a}
\end{eqnarray*}
for $u\in V$ and $w_{a}\in W^{a}$. 
Since $V'$ as a $V$-module is isomorphic to $V$, we have 
$e'=e$. From \cite{FHL}, we know that there is a nondegenerate
invariant 
bilinear form $(\cdot, \cdot)$ on $V$ such that $(\mathbf{1}, 
\mathbf{1})=1$. 
We choose $\Y_{aa'; 1}^{e}=\Y_{aa'; 1}^{e'}$
to be the intertwining operator defined using the action of 
$\sigma_{23}$ by
$$\Y_{aa'; 1}^{e'}=\sigma_{23}(\Y_{ae; 1}^{a}),$$
that is,
$$(u, \Y_{aa'; 1}^{e'}(w_{a}, x)w'_{a})
=e^{\pi i h_{a}}\langle \Y_{ae; 1}^{a}(e^{xL(1)}(e^{-\pi i}x^{-2})^{L(0)}w_{a}, x^{-1})u, 
w'_{a}\rangle$$
for $u\in V$, $w_{a}\in W^{a}$ and $w'_{a'}\in (W^{a})'$. Since the actions of
$\sigma_{12}$
and $\sigma_{23}$ generate the action of $S_{3}$ on $\mathcal{V}$, we have
$$\Y_{a'a; 1}^{e}=\sigma_{12}(\Y_{aa'; 1}^{e})$$
for any $a\in \mathcal{A}$.
When $a_{1}, a_{2}, a_{3}\ne e$, we choose 
$\mathcal{Y}_{a_{1}a_{2}; i}^{a_{3}}$, $i=1, \dots, 
N_{a_{1}a_{2}}^{a_{3}}$, to be an arbitrary basis
of $\mathcal{V}_{a_{1}a_{2}}^{a_{3}}$.
Note that for each element 
$\sigma\in S_{3}$, 
$\{\sigma(\mathcal{Y})_{a_{1}a_{2}; i}^{a_{3}}\;|\;i=1, \dots, 
N_{a_{1}a_{2}}^{a_{3}}\}$ is also 
a basis of $\mathcal{V}_{a_{1}a_{2}}^{a_{3}}$.

We now discuss modular transformations. 
Let $q_{\tau}=e^{2\pi i\tau}$
for $\tau \in \mathbb{H}$ ($\mathbb{H}$ is the upper-half plane). 
We consider the $q_{\tau}$-traces of the vertex operators 
$Y_{W^{a}}$ for $a\in \mathcal{A}$ on 
the irreducible $V$-modules $W^{a}$ of the following form:
\begin{equation}\label{1-trace}
\tr_{W^{a}}
Y_{W^{a}}(e^{2\pi iz L(0)}u, e^{2\pi i z})
q_{\tau}^{L(0)-\frac{c}{24}}
\end{equation}
for $u\in V$. 
In \cite{Z}, under some conditions slightly different 
from (mostly stronger than) those we assume in this paper, 
Zhu proved that these $q$-traces are independent of 
$z$, are absolutely convergent
when $0<|q_{\tau}|<1$ and can be analytically extended to 
analytic functions of $\tau$ in the upper-half plane. 
We shall denote the analytic extension of (\ref{1-trace})
by 
$$E(\tr_{W^{a}}
Y_{W^{a}}(e^{2\pi iz L(0)}u, e^{2\pi i z})
q_{\tau}^{L(0)-\frac{c}{24}}).$$
In \cite{Z}, under his conditions alluded to above,
Zhu also proved the following modular invariance property:
For 
$$\left(\begin{array}{cc}a&b\\ c&d\end{array}\right)\in SL(2, \Z),$$
let $\tau'=\frac{a\tau+b}{c\tau+d}$. Then there exist 
unique $A_{a_{1}}^{a_{2}}\in \C$ for $a_{1}, a_{2}\in \mathcal{A}$
such that 
\begin{eqnarray*}
\lefteqn{E\left(\tr_{W^{a_{1}}}
Y_{W^{a_{1}}}\left(e^{\frac{2\pi iz}{c\tau+d}L(0)}
\left(\frac{1}{c\tau+d}\right)^{L(0)}
u, e^{\frac{2\pi i z}{c\tau+d}}\right)
q_{\tau'}^{L(0)-\frac{c}{24}}\right)}\nn
&&=\sum_{a_{2}\in \A}A_{a_{1}}^{a_{2}}
E(\tr_{W^{a_{2}}}
Y_{W^{a_{2}}}(e^{2\pi iz L(0)}u, e^{2\pi i z})
q_{\tau}^{L(0)-\frac{c}{24}})
\end{eqnarray*}
for $u\in V$. In \cite{DLM}, Dong, Li and Mason, 
among many other things,
improved Zhu's results above by showing that the results 
of Zhu above also hold for vertex operator algebras satisfying 
the conditions (slightly weaker than what) we assume in this paper. 
In particular, for 
$$\left(\begin{array}{cc}0&1\\ -1&0\end{array}\right)\in SL(2, \Z),$$
there exist unique $S_{a_{1}}^{a_{2}}\in \C$ for $a_{1}\in \mathcal{A}$
such that 
\begin{eqnarray*}
\lefteqn{E\left(\tr_{W^{a_{1}}}
Y_{W^{a_{1}}}\left(e^{-\frac{2\pi iz}{\tau}L(0)}
\left(-\frac{1}{\tau}\right)^{L(0)}
u, e^{-\frac{2\pi i z}{\tau}}\right)
q_{-\frac{1}{\tau}}^{L(0)-\frac{c}{24}}\right)}\nn
&&=\sum_{a_{2}\in \A}S_{a_{1}}^{a_{2}}
E(\tr_{W^{a_{2}}}
Y_{W^{a_{2}}}(e^{2\pi iz L(0)}u, e^{2\pi i z})
q_{\tau}^{L(0)-\frac{c}{24}})
\end{eqnarray*}
for $u\in V$. When $u=\mathbf{1}$, we see that the matrix
$S=(S_{a_{1}}^{a_{2}})$ actually acts on the space 
of spanned by the 
vacuum characters $\tr_{W^{a}}q_{\tau}^{L(0)-\frac{c}{24}}$
for $a\in \mathcal{A}$.

\renewcommand{\theequation}{\thesection.\arabic{equation}}
\renewcommand{\thethm}{\thesection.\arabic{thm}}
\setcounter{equation}{0}
\setcounter{thm}{0}

\section{The Verlinde conjecture and consequences}

In \cite{H7}, we proved the following general version of 
the Verlinde conjecture in the framework of vertex operator algebras
(cf. Section 3 in \cite{V} and Section 4 in 
\cite{MS1}):

\begin{thm}\label{main}
Let $V$ be a vertex operator algebra satisfying the 
following conditions:

\begin{enumerate}

\item $V_{(n)}=0$ for $n<0$, $V_{(0)}=\mathbb{C}\mathbf{1}$ and 
$V'$ is isomorphic to $V$ as a $V$-module.

\item Every $\mathbb{N}$-gradable weak $V$-module is completely 
reducible.

\item $V$ is $C_{2}$-cofinite, that is, $\dim V/C_{2}(V)<\infty$. 

\end{enumerate}
Then for $a\in \A$, 
$$F(\Y_{ae; 1}^{a} \otimes \Y_{a'a; 1}^{e};
\Y_{ea; 1}^{a}\otimes \Y_{aa'; 1}^{e})\ne 0$$
and 
\begin{equation}\label{diag2}
\sum_{a_{1}, a_{3}\in \mathcal{A}}(S^{-1})_{a_{4}}^{a_{1}}
N_{a_{1}a_{2}}^{a_{3}}
S_{a_{3}}^{a_{5}}=\delta_{a_{4}}^{a_{5}}
\frac{(B^{(-1)})^{2}(\Y_{a_{4}e; 1}^{a_{4}}\otimes \Y_{a'_{2}a_{2}; 1}^{e};
\Y_{a_{4}e; 1}^{a_{4}}\otimes \Y_{a'_{2}a_{2}; 1}^{e})}
{F(\Y_{a_{2}e; 1}^{a_{2}} \otimes \Y_{a'_{2}a_{2}; 1}^{e};
\Y_{ea_{2}; 1}^{a_{2}}\otimes \Y_{a_{2}a'_{2}; 1}^{e})}.
\end{equation}
In particular, the matrix $S$ diagonalizes the
matrices $\mathcal{N}(a_{2})$ for all  $a_{2}\in \A$.
\end{thm}
{\it Sketch of the proof.}\hspace{2ex}
Moore and Seiberg showed in \cite{MS1}
that the conclusions of the theorem follow
from the following formulas (which they 
derived by assuming the axioms of rational conformal field
theories): For $a_{1}, a_{2}, a_{3}\in \A$, 
\begin{eqnarray}\label{formula1-2}
\lefteqn{\sum_{i=1}^{N_{a_{1}a_{2}}^{a_{3}}}\sum_{k=1}^{N_{a'_{1}a_{3}}^{a_{2}}}
F(\Y_{a_{2}e; 1}^{a_{2}}\otimes \Y_{a'_{3}a_{3}; 1}^{e}; 
\Y_{a'_{1}a_{3}; k}^{a_{2}}\otimes \Y_{a_{2}a'_{3}; i}^{a'_{1}})\cdot}\nn
&&\quad\quad\quad\quad\quad\cdot 
F(\Y_{a'_{1}a_{3}; k}^{a_{2}}\otimes \sigma_{123}(\Y_{a_{2}a'_{3}; i}^{a'_{1}});
\Y_{ea_{2}; 1}^{a_{2}}\otimes \Y_{a'_{1}a_{1}; 1}^{e})\nn
&&=N_{a_{1}a_{2}}^{a_{3}}
F(\Y_{a_{2}e; 1}^{a_{2}} \otimes \Y_{a'_{2}a_{2}; 1}^{e};
\Y_{ea_{2}; 1}^{a_{2}}\otimes \Y_{a_{2}a'_{2}; 1}^{e})
\end{eqnarray}
and
\begin{eqnarray}\label{formula2}
\lefteqn{\sum_{a_{4}\in \mathcal{A}}S_{a_{1}}^{a_{4}}
(B^{(-1)})^{2}(\Y_{a_{4}e; 1}^{a_{4}}\otimes \Y_{a'_{2}a_{2}; 1}^{e};
\Y_{a_{4}e; 1}^{a_{4}}\otimes \Y_{a'_{2}a_{2}; 1}^{e})
(S^{-1})_{a_{4}}^{a_{3}}}\nn
&&=\sum_{i=1}^{N_{a_{1}a_{2}}^{a_{3}}}\sum_{k=1}^{N_{a'_{1}a_{3}}^{a_{2}}}
F(\Y_{a_{2}e; 1}^{a_{2}}\otimes \Y_{a'_{3}a_{3}; 1}^{e}; 
\Y_{a'_{1}a_{3}; k}^{a_{2}}\otimes \Y_{a_{2}a'_{3}; i}^{a'_{1}})\cdot\nn
&&\quad\quad\quad\quad\quad\cdot 
F(\Y_{a'_{1}a_{3}; k}^{a_{2}}\otimes \sigma_{123}(\Y_{a_{2}a'_{3}; i}^{a'_{1}});
\Y_{ea_{2}; 1}^{a_{2}}\otimes \Y_{a'_{1}a_{1}; 1}^{e}).
\end{eqnarray}
So the main work is to prove these two formulas.  The proofs of these
formulas in \cite{H7} are based in turn on the proofs of a number of other 
formulas and on nontrivial applications of a
number of results in the theory of vertex operator algebras, so here we
can only outline what is used in the proofs. 

The proof of the first formula (\ref{formula1-2})
uses mainly the works of Lepowsky and the author \cite{HL1}--\cite{HL4}
and of 
the author \cite{H1} \cite{H2} \cite{H4} and 
\cite{H5} on the tensor product theory,
intertwining operator algebras and the construction 
of genus-zero chiral conformal field theories. 
The main technical result used is the associativity 
for intertwining operators proved in \cite{H1} and \cite{H5}
for vertex operator algebras satisfying the three conditions 
stated in the theorem.
Using the associativity for 
intertwining operators  repeatedly
to express the correlation functions obtained 
from products of 
three suitable intertwining operators as linear 
combinations of the correlation functions obtained 
from iterates of 
three intertwining operators in two ways, we obtain a formula 
for the matrix elements of the fusing isomorphisms. Then using 
certain properties of the matrix elements of the fusing 
isomorphisms and their inverses proved in 
\cite{H7}, we obtain the first formula 
(\ref{formula1-2}).

The proof of the second formula (\ref{formula2}) 
heavily uses the results obtained in 
\cite{H6} on the convergence and analytic extensions 
of the $q_{\tau}$-traces of products of what we call ``geometrically-modified 
intertwining operators'', the genus-one associativity, and the modular
invariance of  these analytic extensions 
of the $q_{\tau}$-traces, where $q_{\tau}=e^{2\pi \tau}$. 
These results allows us to (rigorously) establish  a 
formula which corresponds to the fact that the modular transformation 
$\tau\mapsto -1/\tau$ changes one basic Dehn twist on the Teichm\"{u}ller
space of genus-one Riemann surfaces to 
the other. Calculating the matrices corresponding to the Dehn twists
and substituting the results into this formula, we
obtain (\ref{formula2}).

As in \cite{MS1}, the conclusions of the theorem follow immediately from
(\ref{formula1-2}) and (\ref{formula2}).
\epfv

\begin{rema}
{\rm Note that finitely generated $\N$-gradable weak $V$-modules are 
what naturally appear in the proofs of 
the theorems on genus-zero and genus-one 
correlation functions. Thus Condition 2 is natural and necessary because 
the Verlinde conjecture concerns $V$-modules,
not finitely generated $\N$-gradable weak $V$-modules. 
Condition 3 would be a consequence of the finiteness 
of the dimensions of genus-one conformal blocks, if the conformal 
field theory had been constructed, and  is thus
natural and necessary. For vertex operator algebras 
associated to affine Lie algebras (Wess-Zumino-Novikov-Witten
models) and vertex operator algebras associated to the Virasoro algebra
(minimal models), Condition 2 can be verified easily by reformulating 
the corresponding complete reducibility results in terms of 
the representation theory
of affine Lie algebras and the Virasoro algebra. For these vertex operator 
algebras, Condition 3 can also be easily verified  by using results in the 
representation theory of affine Lie algebras and the Virasoro algebra.
In fact, Condition 3 was stated to hold for these algebras in Zhu's paper
\cite{Z}  and was verified by Dong-Li-Mason 
\cite{DLM} (see also \cite{AN} for the case of minimal models).}
\end{rema}

Using  the fact that 
$N_{ea_{1}}^{a_{2}}=\delta_{a_{1}}^{a_{2}}$ 
for $a_{1}, a_{2}\in \mathcal{A}$, we can easily 
derive the following formulas from Theorem \ref{main} (cf.
Section 3 in \cite{V}):

\begin{thm}
Let $V$ be a vertex operator algebra satisfying the conditions
in Section 1. Then we have $S_{e}^{a}\ne 0$ for $a\in \mathcal{A}$
and
\begin{equation}\label{v-form}
N_{a_{1}a_{2}}^{a_{3}}=
\sum_{a_{4}\in \A}\frac{S_{a_{1}}^{a_{4}}S_{a_{2}}^{a_{4}}S_{a_{4}}^{a'_{3}}}
{S_{0}^{a_{4}}}.
\end{equation}
\end{thm}

\begin{thm}
For $a_{1}, a_{2}\in \mathcal{A}$, 
\begin{equation}\label{s-form-3}
S_{a_{1}}^{a_{2}}
=\frac{S_{e}^{e}((B^{(-1)})^{2}(\Y_{a_{2}e; 1}^{a_{2}}
\otimes \Y_{a'_{1}a_{1}; 1}^{e};
\Y_{a_{2}e; 1}^{a_{2}}\otimes \Y_{a'_{1}a_{1}; 1}^{e}))}
{F(\Y_{a_{1}e; 1}^{a_{1}} \otimes \Y_{a'_{1}a_{1}; 1}^{e};
\Y_{ea_{1}; 1}^{a_{1}}\otimes \Y_{a_{1}a'_{1}; 1}^{e})
F(\Y_{a_{2}e; 1}^{a_{2}} \otimes \Y_{a'_{2}a_{2}; 1}^{e};
\Y_{ea_{2}; 1}^{a_{2}}\otimes \Y_{a_{2}a'_{2}; 1}^{e})}.
\end{equation}
\end{thm}

Using (\ref{s-form-3}) and certain properties of 
the matrix elements of the fusing and braiding isomorphisms
proved in \cite{H7}, 
we can prove the following:

\begin{thm}
The matrix $(S_{a_{1}}^{a_{2}})$ is symmetric. 
\end{thm}

\renewcommand{\theequation}{\thesection.\arabic{equation}}
\renewcommand{\thethm}{\thesection.\arabic{thm}}
\setcounter{equation}{0}
\setcounter{thm}{0}

\section{Rigidity, nondegeneracy property and modular tensor categories}

A tensor category with tensor product bifunctor $\boxtimes$ 
and unit object $V$ 
is rigid if for every object 
$W$ in the category, there are right and left dual objects $W^{*}$
and $^{*}W$ together with morphisms $e_{W}: W^{*}\boxtimes W\to V$,
$i_{W}:V\to W\boxtimes W^{*}$, $e'_{W}: W\boxtimes {}^{*}W\to V$
and $i'_{W}: V\to {}^{*}W\boxtimes W$ 
such that the compositions of the morphisms in the sequence
$$\begin{CD}
W&@>>>
&V\boxtimes W
&@>i_{W}\boxtimes I_{W}>>
&(W\boxtimes W^{*})\boxtimes 
W&@>>>\\
&@>>>&W\boxtimes (W^{*}\boxtimes 
W)&@>I_{W}\boxtimes e_{W}>>&
W\boxtimes V&@>>>&W
\end{CD}$$
and three similar sequences are equal to the identity
$I_{W}$ on $W$. Rigidity is a standard notion in the theory 
of tensor categories. 
A rigid braided tensor category together with a twist (a natural 
isomorphism from the category to itself) satisfying natural conditions 
(see \cite{T} and \cite{BK} for the precise conditions)
is called a ribbon category. A semisimple ribbon category with 
finitely many inequivalent irreducible objects is a modular 
tensor category if it has the following nondegeneracy 
property:
The $m\times m$ matrix formed by the traces
of the morphism $c_{W_{i}W_{j}}\circ c_{W_{j}W_{i}}$ 
in the ribbon category for 
irreducible modules $W_{1}, \dots, W_{m}$ is invertible. 
The term ``modular tensor category'' was first suggested 
by I. Frenkel to summarize Moore-Seiberg's theory of polynomial 
equations. See \cite{T} and \cite{BK} for 
details of the theory of modular tensor 
categories.

The results in the proceeding section give the following:

\begin{thm}
Let $V$ be a simple vertex operator algebra satisfying the 
following conditions:

\begin{enumerate}

\item $V_{(n)}=0$ for $n<0$, $V_{(0)}=\mathbb{C}\mathbf{1}$, 
$W_{(0)}=0$ for any irreducible $V$-module which is 
not equivalent to $V$.

\item Every $\mathbb{N}$-gradable weak $V$-module is completely 
reducible.

\item $V$ is $C_{2}$-cofinite, that is, $\dim V/C_{2}(V)<\infty$. 

\end{enumerate}
Then the braided tensor category structure on the 
category of $V$-modules constructed in 
\cite{HL1}--\cite{HL4}, \cite{H1} and \cite{H5}
is rigid, has a natural structure of 
ribbon category and has the nondegeneracy property.  
In particular, the category of $V$-modules 
has a natural structure of modular tensor category.
\end{thm}
{\it Sketch of the proof.}\hspace{2ex}
Note that Condition 1 implies that
$V'$ is equivalent to $V$ as a $V$-module. Thus Condition 1 
is stronger than Condition 1 in the preceding section. In particular,
we can use
all the results in the proceeding section. This slightly stronger
Condition 1 is needed in the proof of the rigidity and 
nondegeneracy property. 
 
We take both the left and right dual of a $V$-module $W$
to be the contragredient module $W'$ of $W$. Since our 
tensor category is semisimple, to prove the rigidity, we need 
only discuss irreducible modules. For any $V$-module
$W=\coprod_{n\in \mathbb{Q}}W_{(n)}$, we use 
$\overline{W}$ to denote its algebraic completion 
$\prod_{n\in \mathbb{Q}}W_{(n)}$. For $a\in \mathcal{A}$,
using the universal property (see Definition 3.1 in 
\cite{HL3} and Definition 12.1 in \cite{HL4})
for the tensor product module $(W^{a})'\boxtimes W^{a}$,
we know that there exists a unique module map $\hat{e}_{a}:
(W^{a})'\boxtimes W^{a}\to V$ such that 
$$\overline{\hat{e}}_{a}(w'_{a}\boxtimes w_{a})
=\Y_{a'a; 1}^{e}(w'_{a}, 1)w_{a}$$
for $w_{a}\in W^{a}$ and $w'_{a}\in (W^{a})'$, where 
$w'_{a}\boxtimes w_{a}\in \overline{(W^{a})'\boxtimes W^{a}}$
is the tensor product of $w_{1}$ and $w_{2}$, 
$\overline{\hat{e}}_{a}: \overline{(W^{a})'\boxtimes W^{a}}\to 
\overline{V}$ is the natural extension of $\hat{e}_{a}$ to 
$\overline{(W^{a})'\boxtimes W^{a}}$. 
Similarly, we have a module map from
$W^{a}\boxtimes (W^{a})'$ to $V$. Since 
$W^{a}\boxtimes (W^{a})'$ is completely reducible and 
the fusion rule $N_{W^{a}(W^{a})'}^{V}$ is $1$, 
there is a $V$-submodule of $W^{a}\boxtimes (W^{a})'$
which is isomorphic to $V$ under the module map from 
$W^{a}\boxtimes (W^{a})'$ to $V$. 
Thus we obtain a module map
$i_{a}: V\to W^{a}\boxtimes (W^{a})'$ which maps
$V$ isomorphically to this submodule of $W^{a}\boxtimes (W^{a})'$. 
Now
\begin{equation}\label{rigid-map}
\begin{CD}
W^{a}&@>>>
&V\boxtimes W^{a}
&@>i_{a}\boxtimes I_{W^{a}}>>
&(W^{a}\boxtimes (W^{a})')\boxtimes 
W^{a}&@>>>\\
&@>>>&W^{a}\boxtimes ((W^{a})'\boxtimes 
W^{a})&@>I_{W^{a}}\boxtimes \hat{e}_{a}>>&
W^{a}\boxtimes V&@>>>&W_{a}
\end{CD}
\end{equation}
is a module map from an irreducible module to itself.
So it must be the identity map multiplied by a number.
One can calculate this number explicitly and it is equal to 
$$F(\Y_{ae; 1}^{a} \otimes \Y_{a'a; 1}^{e};
\Y_{ea; 1}^{a}\otimes \Y_{aa'; 1}^{e}).$$
From Theorem \ref{main}, this number is not $0$. 
Let 
$$e_{a}=\frac{1}{F(\Y_{ae; 1}^{a} \otimes \Y_{a'a; 1}^{e};
\Y_{ea; 1}^{a}\otimes \Y_{aa'; 1}^{e})}\hat{e}_{a}$$
Then the map obtained from (\ref{rigid-map}) by 
replacing $\hat{e}_{a}$ by $e_{a}$ 
is the identity. 
Similarly, we can prove that all the other maps 
in the definition of rigidity are also equal to the 
identity. Thus the tensor category is rigid.

For any $a\in \A$, we define the twist on $W^{a}$ 
to $e^{2\pi i h_{a}}$. Then it is easy to verify that 
the rigid braided tensor category with this twist is 
a ribbon category.

To prove the nondegeneracy property, we use the formula 
(\ref{s-form-3}). Now it is easy to calculate 
in the tensor category
the trace of $c_{W^{a_{2}}, W^{a_{1}}}\circ 
c_{W^{a_{1}}, W^{a_{2}}}$ for $a_{1}, a_{2}\in \mathcal{A}$, 
where $c_{W^{a_{1}}, W^{a_{2}}}: W^{a_{1}}\boxtimes 
W^{a_{2}}\to W^{a_{2}}\boxtimes 
W^{a_{1}}$ is the braiding 
isomorphism. The result is 
$$\frac{((B^{(-1)})^{2}(\Y_{a_{2}e; 1}^{a_{2}}
\otimes \Y_{a'_{1}a_{1}; 1}^{e};
\Y_{a_{2}e; 1}^{a_{2}}\otimes \Y_{a'_{1}a_{1}; 1}^{e}))}
{F(\Y_{a_{1}e; 1}^{a_{1}} \otimes \Y_{a'_{1}a_{1}; 1}^{e};
\Y_{ea_{1}; 1}^{a_{1}}\otimes \Y_{a_{1}a'_{1}; 1}^{e})
F(\Y_{a_{2}e; 1}^{a_{2}} \otimes \Y_{a'_{2}a_{2}; 1}^{e};
\Y_{ea_{2}; 1}^{a_{2}}\otimes \Y_{a_{2}a'_{2}; 1}^{e})}.
$$
By (\ref{s-form-3}), this is equal to 
$$\frac{S_{a_{1}}^{a_{2}}}{S_{e}^{e}},$$
and these numbers form an invertible matrix. 
The other data and axioms for modular tensor 
categories can be given or proved trivially.
Thus the tensor category is 
modular.  The details 
will be given in  \cite{H8}.
\epfv

\noindent {\small \sc Department of Mathematics, Rutgers University,
110 Frelinghuysen Rd., Piscataway, NJ 08854-8019}

\noindent {\em E-mail address}: yzhuang@math.rutgers.edu

\end{document}